\documentclass[a4paper,12pt]{article}

\usepackage{amsmath}
\usepackage{amsthm}
\usepackage{amssymb}
\usepackage{latexsym}

\newtheorem{thm}{Theorem}[section]
\newtheorem{cor}[thm]{Corollary}
\newtheorem{pro}[thm]{Proposition}
\newtheorem{lem}[thm]{Lemma}
\newtheorem{rem}[thm]{Remark}

\newcommand{\mbf}{\mathbf}
\newcommand{\mrm}{\mathrm}
\newcommand{\mcl}{\mathcal}
\newcommand{\mbb}{\mathbb}

\newcommand{\bfs}{\bfseries}
\newcommand{\non}{\nonumber}

\newcommand{\dis}{\displaystyle}

\newcommand{\lf}{\lfloor}
\newcommand{\rf}{\rfloor}
\newcommand{\bsl}{\backslash}
\newcommand{\emp}{\emptyset}
\newcommand{\sub}{\subset}
\newcommand{\bou}{\partial}
\newcommand{\lra}{\longleftrightarrow}

\newcommand{\wt}{\widetilde}
\newcommand{\lc}{\lceil}
\newcommand{\rc}{\rceil}

\newcommand{\Lam}{\Lambda}
\newcommand{\lam}{\lambda}
\newcommand{\del}{\delta}
\newcommand{\Del}{\Delta}
\newcommand{\Om}{\Omega}
\newcommand{\om}{\omega}
\newcommand{\si}{\sigma}
\newcommand{\ga}{\gamma}
\newcommand{\Ga}{\Gamma}
\newcommand{\al}{\alpha}
\newcommand{\var}{\varepsilon}

\newcommand{\Z}{\mbb{Z}}
\newcommand{\R}{\mbb{R}}
\newcommand{\N}{\mbb{N}}
\newcommand{\F}{\mcl{F}}
\newcommand{\B}{\mbb{B}}
\newcommand{\E}{\mbb{E}}
\newcommand{\p}{\mbb{P}}

\newcommand{\mclp}{\mcl{P}}

\newcommand{\CP}{\mrm{CP}}
\newcommand{\ind}{\mbf{1}}

\newcommand{\half}{\frac{1}{2}}
\newcommand{\cov}{\mrm{cov}}

\newcommand{\dlra}{\stackrel{*}{\lra}}

\title{On Traversable Length inside Semi-Cylinder\\ 
in 2D supercritical Bond Percolation}
\author{
Nobuaki Sugimine\thanks{
Corresponding author. Tel.:079-565-8300; Fax.:079-565-9077; 
E-mail: sugimine@ksc.kwansei.ac.jp} \\ 
\small{Department of Physics, Faculty of Science and Technology,} \\
\small{Kwansei Gakuin University,}\\
\small{2-1 Gakuen Sanda, 669-1337, Japan}\\
and\\
Masato Takei\thanks{E-mail: takei@math.keio.ac.jp} \\
\small{Department of Mathematics, 
Faculty of Science and Technology,} \\
\small{Keio University,}\\
\small{3-14-1 Hiyoshi Kohoku-ku Yokohama, 223-8552, Japan}
}
\date{}
\begin{document}
\baselineskip=7mm
\maketitle
\begin{center}{\bfs Abstract}
\end{center}
{\small{
We investigate a limit theorem on traversable length inside semi-cylinder 
in the 2-dimensional supercritical Bernoulli bond percolation, 
which gives an extension of Theorem 2 in \cite{G81}. 
This type of limit theorems 
was originally studied for the extinction time 
for the 1-dimensional contact process on a finite interval in \cite{WA05}. 
Actually, our main result Theorem \ref{thm1} is stated 
under a rather general 2-dimensional bond percolation setting. 
}}
\section{Introduction}
Grimmett\cite{G81} proved that traversable length by open paths 
has logarithmic scale 
in the 2-dimensional subcritical Bernoulli bond percolaiton. 
By the self-duality, this assertion is equivalent to 
that exponential scale length is traversable 
by open paths in the 2-dimensional supercritical Bernoulli bond percolation. 
More precisely, the supercritical version of the assertion is 
the following limit theorem: 
For $p>1/2$, 
\begin{align}
\lim_{N\to\infty}\mclp_p\left(\begin{array}{c}
\mbox{there exists some crossing open path}\\
\mbox{from the bottom to the top in }T(e^{aN},N)
\end{array}\right)
=
\left\{\begin{array}{cl}
1 & \mbox{if }a<\al(1-p),\\
0 & \mbox{if }a>\al(1-p),
\end{array}\right.
\label{intro1}
\end{align}
where $T(M,N)=\{(x_1,x_2)\in\Z^2:1\le x_1\le N,\ 0\le x_2\le M\}$ 
for $M,N\in\N$ and 
\begin{align}
\al(r)=-\lim_{N\to\infty}\frac{1}{N}\log
\mclp_r(
\mbox{$(0,0)$ is connected to the vertical line $x_1=N$ by open paths})
\label{intro2}
\end{align}
for $r\in(0,1]$. 
A similar result as in the subcritical case 
was obtained by Higuchi\cite{H87} 
for a class of site percolation in (strongly) mixing random fields 
on the $d$-dimensional lattice. 
One of typical examples is the 2-dimensional Ising percolation 
in high temperature phase without external magnetic fields.

In the 1-dimensional contact process, 
same types of limit theorems 
for the extinction time of the process on a finite interval 
were proved 
in subcritical region (see \cite{DL88}) and supercritical region 
(see \cite{DS88}). 
A planar graph duality in the graphical representation 
for the contact process plays a central role in \cite{DS88}. 
For this reason, 
the method in \cite{DS88} was applied to the Bernoulli bond percolaiton 
to give another proof of (\ref{intro1}). 
Further, Durrett and Schonmann\cite{DS88} obtained that 
the traversable length $\si_N/\mclp_p[\si_N]$ scaled by its mean 
converges to a mean one exponential distribution 
in the sense of weak convergence, 
where 
\begin{align*}
\si_N=\sup\{ x_2\in\N: 
\mbox{there exists some open path from $I_N$ to $(x_1,x_2)$ in 
$[1,N]\times[0,\infty)$} \}
\end{align*}
and $I_N=\{1,\ldots,N\}\times\{0\}$. 
As for higher dimensional versions of these types of limit theorems, 
we refer to Part I of \cite{Lig99}.
Chen, Liu, and Zhang\cite{CLZ06} carried out similar analysis 
of reversible nearest neighbor particle systems.

On the basis of the argument in \cite{DS88}, 
Wagner and Anantharam\cite{WA05} studied the extinction time 
$\si_N^{\CP}$ 
for the 1-dimensional contact process with piecewise homogeneous birth rates 
and an identical death rate on a finite interval. 
The precise definition is as follows: 
Let the death rates for all vertices be identically equal to 
the normalized rate $1$. 
Divide the interval $[1,N]$ into $K$ intervals $I_{N,i}$'s 
with length $k_iN$'s. 
For every interval $I_{N,i}$, 
the birth rates for all vertices in $I_{N,i}\cap\Z$ are assumed to 
be equal to $\lam_i$. 
One of results in \cite{WA05} affirms that 
if all $\lam_i$'s are larger than the critical point $\lam_c$ 
of the (original) 1-dimensional contact process, 
\begin{align}
\lim_{N\to\infty}\p_N^{\CP}\left( 
\left|\frac{ \log\si_N^{\CP} }{N}
-\sum_{i=1}^Kk_i\ga^{\CP}({\lam_i})\right|>\del 
\right)=0
\label{intro3}
\end{align}
holds for any $\del>0$, where 
\begin{align*}
\ga^{\CP}(\lam)=-\lim_{L\to\infty}\frac{1}{L}\log\mclp_\lam^{\CP}\left( 
\begin{array}{c}
\mbox{the 1-dimensional contact process with $L$ initial}\\
\mbox{particles on $\{1,\ldots,L\}$ eventually extincts}
\end{array} \right).
\end{align*}

In this paper, we consider a similar type of limit theorems as (\ref{intro3}) 
for a class of 2-dimensional bond percolation models 
with the exponential decay of dual connectivity (DC) 
and the ratio weak mixing (RWM). 
Especially, (DC) is a more important notion 
since dual models of 2-dimensional bond percolation 
are also 2-dimensional bond percolation. 
The self-duality holds for the (infinite volume) random-cluster model 
with parameters $(q,p)$ 
(which is the Bernoulli bond percolation when $q=1$) 
in 2-dimensions. 
Using this property, (DC) is proved for the random-cluster model 
with $q=1,2$ and large enough $q$ 
in whole subcritical region. 
It is also believed that the random-cluster model with $q\ge 1$ 
has (DC) in whole subcritical region. 
In addition, for the random-cluster model with $q\ge 1$, 
(RWM) also follows from (DC) 
(see Theorem 3.4 and Remark 3.5 in \cite{Al98}). 

The rest of this paper is structured as follows. 
Our results and some definition are described in Section 2. 
Section 3 is devoted to the proof of Theorem \ref{thm1} 
which goes along the same way as in \cite{WA05} 
except for using (DC) and (RWM) instead of the independence property. 
Section 4 is devoted to the proofs of Theorems 
\ref{thm5}, \ref{thm3}, and \ref{thm4}. 
\section{Results}
\subsection{Main result} 
Let $\E_\Lam=\{\{x,y\}:x,y\in\Lam\mbox{ such that }|x-y|_1=1\}$ 
for every $\Lam\sub\Z^2$ and $\E=\E_{\Z^2}$, 
where $|x-y|_1$ means the $l_1$-distance between $x$ and $y$. 
We take as state space the set $\{0,1\}^{\E}$ 
and denote $(\om_b)_{b\in\E}\in\{0,1\}^{\E}$ by $\om$. 
For $\om\in\{0,1\}^{\E}$,
we declare a bond $b\in\E$ 
to be {\it open} (resp. {\it closed}) ({\it in $\om$}) 
if $\om_b=1$ (resp. $\om_b=0$). 
Let $(b_i)$ be a finite or an infinite sequence of bonds such that 
$b_i\neq b_j$ if $i\neq j$. 
We call such $(b_i)$ a {\it path} 
if $b_1\cap b_2\neq\emp$ and 
$(b_i\bsl b_{i-1})\cap b_{i+1}\neq\emp$ for all $i\ge 2$. 
We call a path $(b_i)$ an {\it open path} ({\it in $\om$}) 
if all bonds $b_i$'s are open. 
For $x,y\in\Z^2$ and $n\in\N$, 
we call a path $(b_i)_{i=1}^n$ an {\it open path from $x$ to $y$} if 
$(b_i)_{i=1}^n$ is an open path such that 
$x\in b_1\bsl b_{2}$ and $y\in b_n\bsl b_{n-1}$. 
For $\Del,\Lam\sub\Z^2$ and $n\in\N$, 
we call a path $(b_i)_{i=1}^n$ an {\it open path from $\Del$ to $\Lam$} if 
$(b_i)_{i=1}^n$ is an open path from $x$ to $y$ 
for some $x\in\Del$ and $y\in\Lam$.
We denote by $\{ \Del\lra\Lam \}$ the event where such a path exists. 
For $\Del\sub\Z^2$, we define $\{ \Del\lra\infty \}$ as the event that 
there exists some infinite open path $(b_i)_{i=1}^\infty$ 
with $x\in b_1\bsl b_{2}$ for some $x\in\Del$. 
In notation below, we often replace $\{x\}$ with $x$. 
For $\Lam\sub\R^2$, we call a path $(b_i)$ a {\it path in $\Lam$} 
when $b_i\in\E_{\Lam\cap\Z^2}$ for all $i$. 
Similarly, we add '{\it in $\Lam$}' to the other terminology above. 
Let $b^*$ denote the dual bond of $b\in\E$. 
We declare the dual bond $b^*$ to be open if and only if $b$ is closed. 
We define a {\it dual open path} and some related notions 
in a similar way as in the case of an open path. 
For sets $\Del$ and $\Lam$ of the dual lattice $(\Z^2)^*$, 
we denote by $\{ \Del\dlra\Lam \}$ the event that 
there exists some dual open path from $x$ to $y$ 
for some $x\in\Del$ and $y\in\Lam$. 

Let $\Phi$ be a 2-dimensional bond percolation model, 
which is a probability measure on $\{0,1\}^{\E}$. 
For every $\Lam\sub\Z^2$, let $\F_\Lam$ denote the $\si$-field 
generated by $\{\om_b:b\in\E_\Lam\}$. 
We say that $\Phi$ possesses the 
{\it{bounded energy property}} (BE) if 
there exists some $r\in(0,1)$ such that for any $b\in\E$, 
$$r\le\Phi(\om_b=1\mid \F_{b^c})\le 1-r.$$ 
We say that $\Phi$ satisfies the 
{\it{exponential decay of dual connectivity property}} (DC) if 
for some $\zeta,C\in(0,\infty)$ and any $x,y\in(\Z^2)^*$, 
$$\Phi(x\dlra y)\le Ce^{-\zeta|x-y|_1}.$$ 
We say that $\Phi$ satisfies the 
{\it{ratio weak mixing property}} (RWM) if 
there exist some $c,C\in(0,\infty)$ such that 
for any $\Del,\Lam\sub\Z^2$ with $\Del\cap\Lam=\emp$, 
$$\sup\left\{\left|\frac{\Phi(A\cap B)}{\Phi(A)\Phi(B)}-1\right|: 
A\in\F_\Del,\ B\in\F_\Lam,\mbox{ and }\Phi(A)\Phi(B)>0\right\}
\le C\sum_{x\in\Del,\, y\in\Lam}e^{-c|x-y|_1}.$$ 

For $r\in\R$, 
we denote by $\lc r\rc$ and $\lf r\rf$ 
the smallest integer larger than $r$ and 
the largest integer smaller than or equal to $r$, respectively. 
Let $\R^+=[0,\infty)$ and $I_N=\{1,\ldots,N\}\times\{0\}$. 
We consider 
\begin{align}
\si_N=\sup\{ n\in\N: 
I_N\lra\{ 1,\ldots,N \}\times\{n\} \mbox{ in }[1,N]\times\R_+ \}
\label{sigman}
\end{align}
in the following bond percolation model $\p_N$: 
For a given $K\in\N$, let $k_1,\ldots,k_K>0$ with $k_1+\cdots+k_K=1$. 
Define $l_0=0$ and $l_i=k_1+\cdots+k_i$ for every $1\le i\le K$. 
\begin{itemize}
	\item[(P1)] For every $1\le i\le K$, 
$$\p_N(\cdot\mid \F_{\Lam^c})=\Phi_i(\cdot\mid \F_{\Lam^c}) 
\quad\p_N\mbox{-a.s.}$$
for any $\Lam\sub
\bigl(\bigl[\lc l_{i-1}N\rc,\lf l_{i}N\rf \bigr]\times\R^+\bigr)\cap\Z^2$. 
	\item[(P2)] $\p_N$ satisfies the FKG inequality. 
	\item[(P3)] $\p_N$ satisfies (DC). 
\end{itemize}
Here, for every $1\le i\le K$, bond percolation model $\Phi_i$ is assumed to 
possess the translation invariance, the FKG inequality, (BE), (DC) and (RWM). 
Notice that for a fixed $K\in\N$, the constants in (BE), (DC), and (RWM) 
for $\Phi_i$'s can be uniformly chosen, respectively. 

Define 
\begin{align}
\ga_i=-\lim_{N\to\infty}\frac{1}{N}\log\Phi_i\left(
\{ I_N\lra\infty \mbox{ in }\R\times\R_+ \}^c\right).
\label{gammap}
\end{align}
Existence of the above limit follows from 
the subadditive argument together with the FKG inequality. 
Further, by (BE) and (DC),
\begin{align*}
\ga_i\in(0,\infty).
\end{align*}
\begin{thm}\label{thm1} 
For any $k_1,\ldots,k_K>0$ with $k_1+\cdots+k_K=1$ and 
$\del>0$, 
\begin{align}
\p_N\left(
\left|\frac{\log\si_N}{N}-\sum_{i=1}^K k_i\ga_i\right|>\del
\right)\to 0
\label{main}
\end{align}
as $N$ goes to infinity. 
\end{thm}
\subsection{Independent bond percolation}
A probability measure $\p$ on $\{0,1\}^{\E}$ is said to be 
{\it independent bond percolation} if every bond becomes open 
independently of all the other bonds. 
Theorem \ref{thm1} immediately leads the following corollary: 
\begin{cor}\label{cor1} 
Consider the Bernoulli bond percolation $\mclp_p$ for $p\in(1/2,1)$. 
Suppose that $\Phi_i=\mclp_{p_i}$ with 
$p_i>1/2$ for every $1\le i\le K$. 
Then, for any $k_1,\ldots,k_K>0$ with $k_1+\cdots+k_K=1$ and 
$\del>0$, (\ref{main}) holds 
for independent bond percolation $\p_N$'s with (P1). 
\end{cor}
\begin{rem}\label{rem1} 
Let \begin{align*}
\ga(p)=-\lim_{N\to\infty}\frac{1}{N}\log\mclp_p\left(
\{ I_N\lra\infty\mbox{ in }\R\times\R_+ \}^c\right)
\end{align*} 
for $p\in[0,1)$. 
Comparing the case where $K=1$ in Corollary \ref{cor1} 
with (\ref{intro1}) (obtained by Grimmett in \cite{G81}),  
we can see that $\ga(p)=\al(1-p)$ 
as Durrett and Schonmann\cite{DS88} pointed out.
Here, $\al(\cdot)$ is the function in (\ref{intro2}). 
According to Theorems 6.10 and 6.14 in \cite{G99}, 
$\al(\cdot)$ is continuous on $(0,1]$.
\end{rem}
Consider a sequence $\{K_N\}_{N\in\N}$ of positive integers 
and let $l_i^{(N)}=i/K_N$ for every $1\le i\le K_N$. 
Define $\{\mrm{Cyl}_N(i)\}_{i=1}^{K_N}$ as follows: 
$\mrm{Cyl}_N(1)=\bigl(-\infty, \lf l_1^{(N)}N\rf\bigr]\times\R$, 
$\mrm{Cyl}_N(K_N)=\bigl[\lc l_{K-1}^{(N)}N\rc,\infty\bigr)\times\R$, 
and $\mrm{Cyl}_N(i)
=\bigl[ \lc l_{i-1}^{(N)}N\rc,\lf l_{i}^{(N)}N\rf \bigr]\times\R$ 
for every $2\le i\le K_N-1$. 
For $b\in\E$ with $b=\{ (x_1,x_2),(y_1,y_2) \}$, 
let $X(b)=\min\{x_1,y_1\}$. 
\begin{thm}\label{thm5}
Consider the Bernoulli bond percolation $\mclp_p$ for $p\in(1/2,1)$ 
and a continuous function $\rho:[0,1]\to(1/2,1)$. 
Take a sequence $\{K_N\}_{N\in\N}$ of positive integers such that 
as $N$ goes to infinity, $K_N,L_N\to\infty$, and 
$L_N^m/N\to\infty$ for some integer $m\ge 2$, 
where $L_N=N/K_N$. 
Let $p_i^{(N)}=\rho(l_i^{(N)})$ for every $1\le i\le K_N$. 
Define $\p_N$ as independent bond percolation such that 
density of edge $b$ is $p_i^{(N)}$ if $X(b)\in\mrm{Cyl}_N(i)$. 
Then, for any $\del>0$, 
\begin{align}
\p_N\left(
\left|\frac{\log\si_N}{N}-\int_0^1\ga(\rho(u))du\right|>\del
\right)\to 0
\label{thm50}
\end{align}
as $N$ goes to infinity. 
Here, $\ga(\cdot)$ is the function in Remark \ref{rem1}. 
\end{thm}
\subsection{Random-cluster models}
Let $q\ge 1$ throughout this paper. 
Let $\om,\xi\in\{0,1\}^{\E}$ and $\Lam\sub\Z^2$. 
A connected component of the graph $(\Z^2,\{ b\in\E:\om_b=1 \})$ 
is called a {\it cluster} ({\it in $\om$}). 
The number of clusters intersecting $\Lam$ is denoted by $k(\om,\Lam)$. 
Let $\om_\Lam\xi$ denote the bond configuration 
such that $(\om_\Lam\xi)_b=\om_b$ if $b\in\E_\Lam$ and 
$(\om_\Lam\xi)_b=\xi_b$ otherwise. 
For a finite set $\Lam\sub\Z^2$ and $p\in[0,1]$, 
the finite volume random-cluster measure $\Phi_{\Lam,p,q}^\xi$ 
on $\{0,1\}^{\E_\Lam}$ 
with the boundary condition $\xi$ is given by 
\begin{align*}
\Phi_{\Lam,p,q}^\xi(\om)=\frac{1}{Z_\Lam^\xi(p,q)}
\left(\prod_{b\in\E_\Lam}p^{\om_b}(1-p)^{1-\om_b}\right)
q^{k(\om_\Lam\xi,\Lam)}
\qquad\mbox{for}\quad\om\in\{0,1\}^{\E_\Lam},
\end{align*}
where $Z_\Lam^\xi(p,q)$ is the normalizing constant.

Taking the thermodynamic limit, 
there exist the infinite volume random-cluster measures 
$\Phi_{p,q}^w$ and $\Phi_{p,q}^f$ 
corresponding to the wired boundary condition $\xi\equiv 1$ 
and the free one $\xi\equiv 0$, respectively. 
The percolation threshold $p_c(q)$ is defined by 
\begin{align*}
p_c(q)
&=\inf\{p\in[0,1]:\Phi_{p,q}^w(O\lra\infty)>0\}\\
&=\inf\{p\in[0,1]:\Phi_{p,q}^f(O\lra\infty)>0\},
\end{align*}
where $O$ indicates the origin of $\Z^2$ 
(see Sections 4 and 5 in \cite{G06}). 
\begin{rem}\label{rem2} 
(i) 
Let $p_{sd}(q)=\sqrt{q}/(1+\sqrt{q})$. 
It holds that $p_c(q)\ge p_{sd}(q)$ and 
there exists a unique infinite volume random-cluster measure $\Phi_{p,q}$ 
for $p\neq p_{sd}(q)$. 
Further,  
$p_c(q)=p_{sd}(q)$ when $q=1,2$ and $q\ge 25.72$ 
(see Sections 6.2 and 6.4 in \cite{G06}).\\ 
(ii) 
If $q=1,2$ or $q\ge 25.72$ and $p>p_c(q)$, 
(DC) holds for the infinite volume random-cluster measure $\Phi_{p,q}$. 
For sufficiently large $p>p_c(q)$, 
(DC) holds for the infinite volume random-cluster measure $\Phi_{p,q}$ 
(see Sections 6.2 and 6.4 in \cite{G06}).\\
(iii) 
In the infinite volume random-cluster measure $\Phi_{p,q}$ 
with $p\neq p_{sd}(q)$, (DC) implies (RWM) 
(see Theorem 3.4 and Remark 3.5 in \cite{Al98}).
\end{rem}

Recall that $X(b)=\min\{x_1,y_1\}$ 
for $b\in\E$ with $b=\{ (x_1,x_2),(y_1,y_2) \}$. 
Let us fix $K\in\N$ and $k_1,\ldots,k_K>0$ with $k_1+\cdots+k_K=1$. 
For $p_1,\ldots,p_K\in[0,1]$ and $N\in\N$, 
let $\mcl{R}_N(p_1,\ldots,p_K;q)$ denote 
the set of all infinite volume random-cluster measures 
defined by the DLR equation which possess 
a cluster-weight $q$ and an edge-weight $p_i$ for every edge $b$ 
with $X(b)\in \mrm{Cyl}(i)$, 
where $\mrm{Cyl}(1)=\bigl(-\infty, \lf l_1N\rf\bigr]\times\R$, 
$\mrm{Cyl}(K)=\bigl[\lc l_{K-1}N\rc,\infty\bigr)\times\R$, and 
$\mrm{Cyl}(i)=\bigl[\lc l_{i-1}N\rc,\lf l_{i}N\rf \bigr]\times\R$ 
for every $2\le i\le K-1$. 
\begin{thm}\label{thm3} 
Consider the infinite volume random-cluster measure $\Phi_{p,q}$ for 
$p\in(p_c(q),1)$. 
For every $1\le i\le K$, suppose that $p_i>p_c(q)$ and 
$\Phi_i=\Phi_{p_i,q}$ satisfies (DC). 
Then, the set $\mcl{R}_N(p_1,\ldots,p_K;q)$ is nonempty 
for any $k_1,\ldots,k_K>0$ with $k_1+\cdots+k_K=1$ and $N\in\N$. 
Moreover, 
(\ref{main}) holds for any $\del>0$ 
if $\p_N\in\mcl{R}_N(p_1,\ldots,p_K;q)$ for all $N\in\N$. 
\end{thm}
\begin{thm}\label{thm4}
Consider the infinite volume random-cluster measure $\Phi_{p,q}$ for 
$p\in(p_c(q),1)$. 
Suppose that $p_i>p_c(q)$ and $\Phi_{p_i,q}$ satisfies (DC) 
for every $1\le i\le K$. 
Consider the semi-cylindrical random-cluster measure $\p_{N,cyl}^w$ 
corresponding to the wired boundary condition 
such that its cluster-weight is $q$ and for every $1\le i\le K$, 
its edge-weight is $p_i$ for every edge $b$ 
with $X(b)\in\bigl[ \lc l_{i-1}N\rc,\lf l_{i}N\rf \bigr]\times\R^+$. 
Then, for any $k_1,\ldots,k_K>0$ with $k_1+\cdots+k_K=1$ and 
$\del>0$, (\ref{main}) holds for $\p_{N,cyl}^w$'s. 
\end{thm}
\section{Proof of Theorem \ref{thm1}}
Although we can prove Theorem \ref{thm1} along the line in \cite{WA05} 
by using (DC) and (RWM) instead of independency, 
we will give its full proof for self-consistency. 
We write $\bar{\ga}=\sum_{i=1}^K k_i\ga_i$. 
We sometimes omit the index $i$ from the notation. 
\subsection{Upper bound}
We will show that 
for any $\del>0$, 
\begin{align}
\lim_{N\to\infty}
\p_N\left(\frac{\log\si_N}{N}>\bar{\ga}+\del\right)=0. 
\label{pro11}
\end{align}
For $M\in\N$, define 
\begin{gather*}
A_N^M=\left\{ I_N\lra\Z\times\{M\}\mbox{ in } \R\times[0,M] \right\}^c,\\
B_N=\left\{ (\half,\half)\dlra(N+\half,\half)\mbox{ in }R_N^+ \right\},
\end{gather*}
and 
\begin{align*}
B_N^M=\left\{ (\half,\half)\dlra(N+\half,\half)\mbox{ in }R_N^+(M) \right\},
\end{align*}
where $R_N^+=[1/2,N+(1/2)]\times[1/2,\infty)$ and 
$R_N^+(M)=[1/2,N+(1/2)]\times[1/2,M-(1/2)]$. 
For every $1\le i\le K$, the following three limits exist 
as in the case of $\ga_i$ (see (\ref{gammap})): 
\begin{gather*}
\ga_i^M=-\lim_{N\to\infty}\frac{1}{N}\log\Phi_i(A_N^M),\\
\mu_i=-\lim_{N\to\infty}\frac{1}{N}\log\Phi_i(B_N),
\end{gather*}
and 
\begin{align*}
\mu_i^M=-\lim_{N\to\infty}\frac{1}{N}\log\Phi_i(B_N^M).
\end{align*}
\begin{lem}\label{lem1}
For every $1\le i\le K$, 
\begin{align*}
\ga_i = \mu_i = \lim_{M\to\infty}\mu_i^M.
\end{align*}
\end{lem}
\begin{proof}
Note that $\ga^M$ and $\mu^M$ are decreasing in $M$. 
By the definitions of $\ga$ and $\ga^M$, 
\begin{align}
\exp(-(\ga+\var)N)
&\le \Phi\left( \{ I_N\lra\infty\mbox{ in }\R\times\R_+ \}^c \right)\non\\
&=\lim_{M\to\infty}\Phi(A_N^M)\non\\
&\le\lim_{M\to\infty}\exp(-\ga^MN)\non
\end{align}
for any $\var>0$ and sufficiently large $N$, 
which together with $\ga\le\ga^M$ implies that 
$\ga=\lim_{M\to\infty}\ga^M$. 
Similarly, $\mu=\lim_{M\to\infty}\mu^M$. 
Further, $\ga^M=\mu^M$ for any $M\in\N$ since 
\begin{align*}
r^2\Phi(B_N^M)\le\Phi(A_N^M)\le r^{-2M}\Phi(B_N^M)
\end{align*}
follow from the FKG inequality and (BE). 
Thus, $\mu=\ga$ holds. 
\end{proof}
\begin{proof}[Proof of the upper bound (\ref{pro11})]
By Lemma \ref{lem1}, there exists some integer $M$ 
such that for all $1\le i\le K$, 
\begin{align}
\mu_i^M\le\ga_i+(\del/6). 
\label{u1}
\end{align}
Take a positive $\eta$ satisfying 
$$\eta\le\min\left\{\frac{\del}{6K\log(1/r)},\ 
\frac{1}{3}\min\{k_1,\ldots,k_K\}\right\}.$$
Let us fix such $\eta$ and $M$. 
For every $1\le i\le K$, consider 
\begin{align*}
B_{N,i}=\left\{ \begin{array}{c}
\bigl( \lc (l_{i-1}+\eta) N\rc-\dis\half,\half \bigr)
\dlra\bigl( \lf (l_{i}-\eta)N\rf-\dis\half,\half \bigr)\\
\mbox{ in }\bigl[\lc l_{i-1} N\rc+\dis\half,\lf l_{i} N\rf-\dis\half\bigr]
\times [\dis\half,M-\dis\half]
\end{array} \right\}
\end{align*}
and 
$$F_N=\{\mbox{all dual bonds in $\B_N$ are open}\},$$
where 
\begin{align*}
\B_N=
&\left\{ \{(j-\half,\half),(j+\half,\half)\}: 
1\le j\le \lf \eta N\rf\mbox{ or }\lc (1-\eta)N\rc\le j\le N \right\}\\
&\cup
\bigcup_{i=1}^{K-1}\left\{\{(j-\half,\half),(j+\half,\half)\}: 
\lc (l_{i}-\eta)N\rc\le j\le \lf (l_{i}+\eta)N\rf\right\}.
\end{align*}
By the FKG inequality, (P1), (RWM) and (\ref{u1}), 
\begin{align}
\p_N\left(B_N^{M}\right)
\ge&
\p_N(F_N) 
\prod_{i=1}^K \p_N\left(B_{N,i}\right)
\ge\half r^{2K\eta N}\prod_{i=1}^K\Phi_i(B_{N,i})
\ge\half
\exp\left(-\left\{\bar{\ga}+(2\del/3)\right\}N\right)
\label{pro12}
\end{align}
for sufficiently large $N$. 
Let $H=\{x\in\Z^2:x_2=iM\mbox{ for some }i\in\N\}$. 
By comparing $\si_N$ with $\si_N$ conditioned by the event 
that all bonds in $\E_H$ are open, 
it is not difficult to see that for any $l\in\N$,  
\begin{align}
\p_N(\si_N>l)
&\le\p_N(\si_N\ge M)^{\lf l/M\rf}
\le\left(1-\p_N\left(B_N^{M}\right)\right)^{\lf l/M\rf}.
\label{pro13}
\end{align}
From (\ref{pro12}) and (\ref{pro13}), we can conclude (\ref{pro11}). 
\end{proof}

\subsection{Lower bound}
Because of (\ref{pro11}), we obtain Theorem \ref{thm1} 
once we can prove that for any $\del>0$, 
\begin{align}
\lim_{N\to\infty}
\p_N\left(\frac{\log\si_N}{N}<\bar{\ga}-\del\right)=0.
\label{pro21}
\end{align}
Let 
\begin{align*}
C_N=\left\{ \mbox{for some $k\in\Z$, }
(\dis\half,\half)\dlra(N+\dis\half,k+\half) \mbox{ in }R_N \right\},
\end{align*}
where $R_N=[1/2,N+(1/2)]\times\R$. 
By Proposition \ref{lem3} mentioned below, 
we can see that for any $l\in\N$ and sufficiently large $N$, 
\begin{align*}
\p_N(\si_N<l)\le (l+1)\p_N(C_N)
\le (l+1)\exp\left(
-\left\{\bar{\ga}-\frac{\del}{2}\right\}N\right),
\end{align*}
which implies (\ref{pro21}). 
\begin{pro}\label{lem3} 
\begin{align}
\lim_{N\to\infty}\frac{1}{N}\log\p_N\left(C_N\right)=-\bar{\ga}.
\label{lem31}
\end{align}
\end{pro}

We prepare some notation and lemmas to prove (\ref{lem31}). 
Let $R_N(M)=[1/2,N+(1/2)]\times[-M+(1/2),M-(1/2)]$ for $M\in\N$. 
Define 
\begin{gather*}
C_N^M=\left\{ 
\mbox{for some }k\in\Z, 
(\dis\half,\half)\dlra(N+\dis\half,k+\half)\mbox{ in }R_N(M) \right\},\\
D_N=\left\{ 
\mbox{for some }j,k\in\Z, 
(\dis\half,j+\half)\dlra(N+\dis\half,k+\half)\mbox{ in }R_N \right\},\\
D_N^M=\left\{ 
\mbox{for some }j,k\in\Z, (\dis\half,j+\half)\dlra(N+\dis\half,k+\half)
\mbox{ in }R_N(M) \right\},\\
E_N=\left\{ 
(\dis\half,\half)\dlra(N+\dis\half,\half)\mbox{ in }R_N \right\}
\end{gather*}
and 
\begin{align*}
E_N^M=\left\{ 
(\dis\half,\half)\dlra(N+\dis\half,\half)\mbox{ in }R_N(M)\right\}.
\end{align*}
\begin{lem}\label{lem4}
For every $1\le i\le K$, 
\begin{align}
\sup_{N\in\N}\frac{1}{N}\log\Phi_i(E_N)
=\lim_{N\to\infty}\frac{1}{N}\log\Phi_i(E_N)
=\lim_{M\to\infty}\lim_{N\to\infty}\frac{1}{N}\log\Phi_i(E_N^M)=-\ga_i. 
\label{lem41}
\end{align}
\end{lem}
\begin{proof}
Note that 
\begin{align*}
B_N^{2M}
&\supset\left\{ 
(\dis\half,M+\half)\dlra(N+\dis\half,M+\half)\mbox{ in }R_N(M)+(0,M) 
\right\}\\
&\qquad
\cap\{\mbox{all dual bonds in $\mbb{B}_N$ are open}\},
\end{align*}
where 
\begin{align*}
\B_N
=&\left\{ \{(\half,\half),(\frac{3}{2},\half)\},
\{(N-\half,\half),(N+\half,\half)\} \right\}\\
&\cup\left\{ \{(\beta,j-\half),(\beta,j+\half)\}: 
1\le j\le M\mbox{ and }\beta=\frac{3}{2},N-\half \right\}.
\end{align*}
Then, by the FKG inequality and (BE), 
\begin{align*}
\lim_{N\to\infty}\frac{1}{N}\log\Phi_i(B_N^{2M})
\ge\lim_{N\to\infty}\frac{1}{N}\log\Phi_i(E_N^M)
\end{align*}
for any $M\in\N$. 
From this and the fact that $B_N^M\sub E_N^M$, 
\begin{align*}
\lim_{N\to\infty}\frac{1}{N}\log\Phi_i(B_N^M)
=\lim_{N\to\infty}\frac{1}{N}\log\Phi_i(E_N^M) 
\end{align*}
for any $M\in\N$. 
Thus, we can obtain (\ref{lem41}) as in the proof of Lemma \ref{lem1}. 
\end{proof}
\begin{lem}\label{lem5}
For every $1\le i\le K$, 
\begin{align}
\lim_{N\to\infty}\frac{1}{N}\log\Phi_i(C_N)=-\ga_i. 
\label{lem51}
\end{align}
\end{lem}
\begin{proof}
By Lemma \ref{lem4} and the fact that $E_N\sub C_N$, 
\begin{align}
-\ga\le\liminf_{N\to\infty}\frac{1}{N}\log\Phi(C_N).
\label{lem52}
\end{align}
By (DC), 
\begin{align}
\Phi(C_N\bsl C_N^M)\le 2CNe^{-\zeta M}
\label{lem53}
\end{align}
for any $M\in\N$. 
From (\ref{lem52}) and (\ref{lem53}), 
\begin{align}
\limsup_{N\to\infty}\frac{1}{N}\log\Phi(C_N)=
\limsup_{N\to\infty}\frac{1}{N}\log\Phi(C_N^{aN})
\label{lem54}
\end{align}
for $a=\lc 6\ga/\zeta\rc$. 
Note that 
\begin{align}
\frac{1}{N}\log\Phi(C_N^{aN})
\le\frac{1}{N}\log(2aN)
+\frac{1}{N}\sup_{|k|\le aN}\log\Phi(C_N(k)), 
\label{lem55}
\end{align}
where 
\begin{align*}
C_N(k)=\left\{ 
(\dis\half,\half)\dlra(N+\half,k+\half)\mbox{ in }R_N \right\}
\end{align*}
for every $k\in\Z$. 
Further, Lemma \ref{lem4} maintains 
\begin{align}
\sup_{k\in\Z}\Phi(C_N(k))\le e^{-\ga N},
\label{lem56}
\end{align}
since by the translation invariance and the FKG inequality, 
$$\Phi(C_N(k))^2\le \Phi(E_{2N})$$
for any $k\in\Z$. 
From (\ref{lem55}) and (\ref{lem56}), 
\begin{align}
\limsup_{N\to\infty}\frac{1}{N}\log\Phi(C_N^{aN})\le -\ga. 
\label{lem57}
\end{align}
Therefore, (\ref{lem51}) follows from (\ref{lem52}), 
(\ref{lem54}), and (\ref{lem57}). 
\end{proof}
\begin{lem}\label{lem6} 
For every $1\le i\le K$, 
\begin{align}
\limsup_{N\to\infty}\frac{1}{N}\log\Phi_i\left(D_N^{N^2}\right)\le -\ga_i.
\label{lem61}
\end{align}
\end{lem}
\begin{proof}
The fact that $\Phi(D_N^{N^2})\le(2N^2+1)\Phi(C_N)$ and Lemma \ref{lem5} 
immediately show (\ref{lem61}). 
\end{proof}
\begin{rem} 
This lemma together with Lemma \ref{lem4} 
means that in (\ref{lem61}), 
the upper limit and the inequality can be replaced with 
limit and equality, respectively. 
\end{rem}
\begin{proof}[Proof of Proposition \ref{lem3}]
From (\ref{pro12}),
\begin{align}
\liminf_{N\to\infty}\frac{1}{N}\log\p_N\left(C_N\right)\ge\bar{\ga}.
\label{lem32}
\end{align}
For every $\eta>0$ as in the proof of (\ref{pro11}) and  
$1\le i\le K$, let 
$$N_i=\lc (l_{i}-\eta) N\rc-\lf (l_{i-1}+\eta)N\rf.$$
Define 
\begin{align*}
D_{N,i}=\left\{\begin{array}{c}
\mbox{for some $j,k\in\Z$, }
\bigl( \lc (l_{i-1}+\eta) N\rc-\dis\half,j+\half \bigr)\dlra
\bigl( \lc (l_{i}-\eta)N\rc-\dis\half,k+\half \bigr)\\
\mbox{in }
\bigl[ \lc (l_{i-1}+\eta) N\rc-\dis\half,\lc (l_{i}-\eta) N\rc-\dis\half \bigr]
\times [0,N_i^2-\dis\half]
\end{array}\right\}.
\end{align*}
By (P1), (RWM), and Lemma \ref{lem6}, 
\begin{align}
\limsup_{N\to\infty}\frac{1}{N}\log\p_N\left(C_N^{aN}\right)
&\le \limsup_{\eta\searrow 0} \limsup_{N\to\infty}\frac{1}{N}
\log\Bigl(2\Phi_1(D_{N,1})\p_N(D_{N,2}\cap\cdots\cap D_{N,K})\Bigr)\non\\
&\le \limsup_{\eta\searrow 0}\sum_{1\le i\le K}
\limsup_{N\to\infty}\frac{1}{N}\log\Phi_i(D_{N,i})\non\\
&=-\bar{\ga} 
\label{lem33} 
\end{align} 
for any $a>0$ and sufficiently large $N$. 
Let $a=\lc 6\ga/\zeta\rc$. 
Then, by (DC) and (\ref{lem32}), 
$$\limsup_{N\to\infty}\frac{1}{N}\log\p_N\left(C_N\right)
=\limsup_{N\to\infty}\frac{1}{N}\log\p_N\left(C_N^{aN}\right),$$
which together with (\ref{lem33}) implies that 
\begin{align}
\limsup_{N\to\infty}\frac{1}{N}\log\p_N\left(C_N\right)
\le -\bar{\ga}. 
\label{lem34}
\end{align}
This and (\ref{lem32}) can lead (\ref{lem31}). 
\end{proof}
\section{Proof of Theorems \ref{thm5}, \ref{thm3}, and \ref{thm4}}
\begin{proof}[Proof of Theorem \ref{thm5}]
Let us fix the integer $m$ as in Theorem \ref{thm5} and an integer $M$. 
Let $p_-=\min\{\rho(u):u\in[0,1]\}$ and $p_+=\max\{\rho(u):u\in[0,1]\}$. 
Note that $1/2<p_-\le p_+<1$. 
We will show that 
\begin{align}
\frac{1}{N}\log\mclp_p(D_N^{N^m})\qquad\mbox{and}\qquad 
\frac{1}{N}\log\mclp_p(B_N^M)
\label{thm51}
\end{align}
are Lipschitz continuous functions in $p$ on $[p_-,p_+]$ 
uniformly in $N\in\N$, 
which implies that 
\begin{align*}
\lim_{N\to\infty}\sum_{i=1}^{K_N}\ga(p_i)k_i=\int_0^1\ga(\rho(u))du
\end{align*}
and both terms in (\ref{thm51}) 
converge uniformly in $p\in[p_-,p_+]$ as $N$ goes to infinity. 
Using these facts, we can obtain (\ref{thm50}) 
in a similar way as in the proof of Theorem \ref{thm1}. 

For simplicity, 
we consider 
\begin{align*}
\frac{1}{N}\log\mcl{P}_p\left(\wt{D}_N\right)\qquad 
\mbox{for}\quad p\in[r_-,r_+], 
\end{align*}
where $r_-=1-p_+$, $r_+=1-p_-$ and 
\begin{align*}
\wt{D}_N=\left\{
\mbox{for some }j,k\in\Z, (0,j)\lra(N,k)\mbox{ in }[0,N]\times[0,2N^m-1]
\right\}.
\end{align*}
Note that $0<r_-\le r_+<1/2$ and $\wt{D}_N$ is a local event. 
Let $\Del_N=([0,N]\times[0,2N^m-1])\cap\Z^2$ and $\E_N=\E_{\Del_N}$. 
By abusing notation, $[0,N]\times[0,2N^m-1]$ is also denoted by $\Del_N$. 
Let $\Om_N$ indicate the number of open bonds in $\E_N$. 
By Russo's formula 
(see Section 2.5 in \cite{G99} or Section 2.4 in \cite{G06}), 
\begin{align}
\frac{d}{dp}\mclp_p\left(\wt{D}_N\right)
=\frac{1}{p(1-p)}\cov_p(\Om_N,\ind_{\wt{D}_N}), 
\label{thm52}
\end{align}
where $\cov_p$ means the covariance with respect to $\mclp_p$ 
and $\ind_{\wt{D}_N}$ denotes the indicator function of $\wt{D}_N$. 

For a set $E\sub\E$, let $|E|$ and $\bou E$ mean 
the cardinality of $E$ and the set of all boundary bonds of $E$, respectively. 
More precisely, $\bou E$ is defined by 
\begin{align*}
\bou E=\{e\in\E: \mbox{$e\notin E$ and $e\cap b\neq\emp$ for some $b\in E$}\}.
\end{align*}
A set $E$ is said to be {\it connected} if for any $b,b'\in E$, 
there exists some path in $E$ which includes both $b$ and $b'$. 
Define the open bond cluster $\wt{C}_{N,x}$ in $\Del_N$ 
(containing $x\in\Z^2$) as follows: 
\begin{align*}
\wt{C}_{N,x}=\{ b\in\E:
\mbox{$b$ is included in some open path from $x$ in $\Del_N$} \}.
\end{align*} 
On the event $\wt{D}_N$, 
there exists some $\wt{C}_{N,(0,j)}$ crossing from the left to the right 
in $\Del_N$. 
Define $\Ga_N$ as $\wt{C}_{N,(0,j)}$ 
with the minimal $j\in\{0,\ldots,2N^m-1\}$ among such $\wt{C}_{N,(0,j)}$'s. 
Then, by the FKG inequality, 
\begin{align}
\cov_p(\Om_N,\ind_{\wt{D}_N})
&=\sum_{C}\mclp_p(\Ga_N=C)\Bigl(
\mclp_p[\Om_N \mid \Ga_N=C]-\mclp_p[\Om_N] \Bigr)\non\\
&\le\sum_{C}\mclp_p(\Ga_N=C)\Bigl((|\E_N|-|C|-|\bou C|)p+|C|-|\E_N|p\Bigr)
\non\\
&\le(1-p)\sum_{n\ge N}\mclp_p(|\Ga_N|\ge n),
\label{thm53}
\end{align}
where $\sum_{C}$ stands for the summation over all connected subsets of 
$\E_N$ crossing from the left to the right in $\Del_N$. 
Note that 
$\mclp_p(|\Ga_N|\ge N)\ge p^N\ge r_-^N$ for all $p\in[r_-,r_+]$.
Let $C_O$ be the open cluster containing the origin $O$ of $\Z^2$. 
In the subcritical regime, the cluster size distribution decays exponentially 
(see Section 6.3 in \cite{G99}). 
This fact together with the FKG inequality implies that 
for some $A\in(0,\infty)$ and all $p\in[r_-,r_+]$, 
\begin{align}
\sum_{n\ge AN}\mclp_p(|\Ga_N|\ge n)
&\le 2N^m\sum_{n\ge AN}\mclp_p(|C_O|\ge \lf n/4 \rf)\non\\
&\le 8N^m\sum_{n\ge AN/4}\mclp_{r_+}(|C_O|\ge n)\non\\
&\le \mclp_p(|\Ga_N|\ge N),
\label{thm54}
\end{align}
where $|C_O|$ means the cardinality of $C_O$. 
From (\ref{thm52})--(\ref{thm54}), 
\begin{align}
\frac{d}{dp}\mclp_p\left(\wt{D}_N\right)
\le\frac{2AN}{r_-}\mclp_p\left(\wt{D}_N\right), 
\label{thm56}
\end{align}
which implies the first term in (\ref{thm51}) 
is uniformly Lipschitz continuous in $N\in\N$. 
As for the second term in (\ref{thm51}), 
the proof is similar as above and easier. 
\end{proof}
\begin{proof}[Proof of Theorem \ref{thm3}]
Let us fix $p_1,\ldots,p_K$ and $q$ as in Theorem \ref{thm3} 
and write $\mcl{R}_N=\mcl{R}_N(p_1,\ldots,p_K;q)$. 
By the definition of $\mcl{R}_N$, 
(P1) holds for any element $\p_N$ of $\mcl{R}_N$. 

The set of all limit random-cluster measures which possess 
a cluster-weight $q$ and an edge-weight $p_i$ 
for every edge $b$ with $X(b)\in \mrm{Cyl}(i)$ 
is denoted by $\mcl{W}_N$. 
The element of $\mcl{W}_N$ corresponding to the wired (resp. free) 
boundary condition 
is denoted by $\p_N^w$ (resp. $\p_N^f$). 
Both measures $\p_N^w$ and $\p_N^f$ satisfy the FKG inequality. 
Further, comparing them with $\Phi_{p_0,q}$ 
in the FKG sense leads their (DC) property, 
where $p_0=\min\{p_1,\ldots,p_K\}$. 
Therefore, there exists a unique infinite cluster almost surely 
under both $\p_N^w$ and $\p_N^f$, 
which implies that  $\p_N^w,\p_N^f\in\mcl{R}_N$ 
(see Section 4.4 in \cite{G06}). 
Thus, $\p_N^w$ and $\p_N^f$ satisfy (P1), (P2), and (P3). 
This together with Theorem \ref{thm1} maintains that 
for a certain $\bar{\ga}$ independent of $w$ and $f$ and any $\del>0$, 
\begin{align*}
\lim_{N\to\infty}\p_N^w\left(
\left|\frac{\log\si_N}{N}-\bar{\ga}\right|>\del\right)
=
\lim_{N\to\infty}\p_N^f\left(
\left|\frac{\log\si_N}{N}-\bar{\ga}\right|>\del\right)
=0. 
\end{align*}
Then, by the FKG inequality, 
\begin{align*}
\lim_{N\to\infty}
\p_N\left(\frac{\log\si_N}{N}>\bar{\ga}+\del\right)
\le 
\lim_{N\to\infty}
\p_N^w\left(\frac{\log\si_N}{N}>\bar{\ga}+\del\right)
=0 
\end{align*}
and
\begin{align*}
\lim_{N\to\infty}
\p_N\left(\frac{\log\si_N}{N}<\bar{\ga}-\del\right)
\le 
\lim_{N\to\infty}
\p_N^f\left(\frac{\log\si_N}{N}<\bar{\ga}-\del\right)
=0. 
\end{align*}
\end{proof}
\begin{proof}[Proof of Theorem \ref{thm4}]
By the FKG inequality, it is sufficient to prove that 
for any $\del>0$, 
\begin{align}
\lim_{N\to\infty}
\p_{N,cyl}^w\left(\frac{\log\si_N}{N}>\bar{\ga}+\del\right)=0. 
\label{thm41}
\end{align}
In the same way as in (\ref{pro12}) and (\ref{pro13}), 
\begin{align*}
\p_{N,cyl}^w(\si_N>l)
&\le\left( 1-\p_{N,cyl}^w\left(B_N^{M}\right) \right)^{\lf l/M\rf}
\le\left( 1-\half\exp\left( -\left\{\bar{\ga}+(2\del/3)\right\}N \right)
\right)^{\lf l/M\rf}
\end{align*}
for some $M\in\N$, which implies (\ref{thm41}). 
\end{proof}

\end{document}